\documentclass[a4paper,11pt,leqno]{amsart}
\usepackage{amsthm}
\usepackage{bm}
\usepackage{amsmath,amssymb,cite,mathrsfs,overpic}

\numberwithin{equation}{section}

\newtheorem{theorem}{Theorem}[section]
\newtheorem{lemma}[theorem]{Lemma}

\theoremstyle{definition}

\newtheorem{example}[theorem]{Example}

\theoremstyle{remark}
\newtheorem{remark}[theorem]{Remark}

\numberwithin{equation}{section}

\allowdisplaybreaks

\long\def\comment#1{}

\begin{document}

\title[Witten multiple zeta-functions]{On Witten multiple zeta-functions associated with semisimple Lie
algebras IV}
\author{Yasushi Komori}
\address{Graduate School of Mathematics, Nagoya University, Chikusa-ku, Nagoya 464-8602 Japan}
\email{komori@math.nagoya-u.ac.jp}

\author{Kohji Matsumoto}
\address{Graduate School of Mathematics, Nagoya University, Chikusa-ku, Nagoya 464-8602 Japan}
\email{kohjimat@math.nagoya-u.ac.jp}

\author{Hirofumi Tsumura}
\address{Department of Mathematics and Information Sciences, Tokyo Metropolitan University, 1-1, Minami-Ohsawa, Hachioji, Tokyo 192-0397 Japan}
\email{tsumura@tmu.ac.jp}

\date{}
\subjclass[2000]{Primary 11M41; Secondary 17B20, 40B05}

\maketitle

\begin{abstract}
In our previous work, we established the theory of multi-variable Witten zeta-functions, which are called the zeta-functions of root systems. We have already considered the cases of types $A_2$, $A_3$, $B_2$, $B_3$ and $C_3$. 
In this paper, we consider the case of $G_2$-type. We define certain analogues of Bernoulli polynomials of $G_2$-type and study the generating functions of them to determine the coefficients of Witten's volume formulas of $G_2$-type. Next we consider the meromorphic continuation of the zeta-function of $G_2$-type and determine its possible singularities. Finally, by using our previous method, we give explicit functional relations for them which include Witten's volume formulas.
\end{abstract}

\baselineskip 19pt

%%%%%%%%%%%%%%%%%%%%%%%%%%%%%%%%%%%%%%
\section{Introduction} \label{sec-1}
%%%%%%%%%%%%%%%%%%%%%%%%%%%%%%%%%%%%%%

Let ${\Bbb N}$ be the set of positive integers, ${\Bbb N}_0={\Bbb N}\cup\{0\}$,
${\Bbb Z}$ the ring of rational integers, ${\Bbb Q}$ the rational number field,
${\Bbb R}$ the real number field, ${\Bbb C}$ the complex number field, 
respectively.

In our previous articles \cite{KMT,KM2,KM3,MTF}, we defined the multi-variable version of Witten zeta-functions, or ``zeta-functions of root systems'', inspired by the original work of Witten \cite{Wi} and of Zagier \cite{Za}. We recall these results as follows.

Let ${\frak g}$ be a complex semisimple Lie algebra with rank $r$, 
$\mathfrak{h}$ be a Cartan subalgebra of $\mathfrak{g}$ and $\mathfrak{h}^*$ be its dual. 
Let $\Delta\subset\mathfrak{h}^*$ be the set of all roots of ${\frak g}$,
% in the vector space $\mathfrak{h}^*$,
%equipped with the inner product $\langle\cdot,\cdot\rangle$
 $\Delta_+$ the set of all positive roots of
${\frak g}$, 
$\Psi=\{\alpha_1,\ldots,\alpha_r\}$ the fundamental
system of $\Delta$, and $\alpha_j^{\vee}$ the coroot associated with
$\alpha_j$ ($1\leq j\leq r$).  Let $\lambda_1,\ldots,\lambda_r$ be the
fundamental weights satisfying
$\lambda_j(\alpha_i^{\vee})
=\delta_{ij}$ (Kronecker's delta).
In the following we
denote the pairing $\lambda(h)$ of $h\in\mathfrak{h}$ and $\lambda\in\mathfrak{h}^*$
by $\langle h,\lambda\rangle$.

In \cite{KMT,KM2}, we defined the multi-variable version of Witten zeta-functions by
\begin{align}
  \label{1-3}
    \zeta_r({\bf s};{\frak g})=\sum_{m_1=1}^{\infty}\cdots
      \sum_{m_r=1}^{\infty}\prod_{\alpha\in\Delta_+}\langle\alpha^{\vee},
      m_1\lambda_1+\cdots+ m_r\lambda_r\rangle ^{-s_{\alpha}}
\end{align}
for ${\bf s}=(s_{\alpha})_{\alpha\in\Delta_+}\in {\Bbb C}^n$, where $n$ is the number of all positive roots. In the case that ${\frak g}$ is of type $X_r$, we call (\ref{1-3}) the zeta-function of the root system of type $X_r$, and also denote it by $\zeta_r({\bf s};X_r)$, where $X=A,B,C,D,E,F,G$. 
Note that the original Witten zeta-function $\zeta_W(s;{\frak g})$, studied by Witten \cite{Wi} and Zagier \cite{Za}, coincides with
\begin{align}
  \label{1-3-2}
    K({\frak g})^s\zeta_r(s,\ldots,s;{\frak g}),
\end{align}
where
\begin{align}
  \label{1-3-3}
   K({\frak g})=\prod_{\alpha\in\Delta_+}\langle\alpha^{\vee},
           \lambda_1+\cdots+\lambda_r\rangle.
\end{align}
Witten's motivation of introducing the above zeta-functions is to 
express the volumes of certain moduli spaces in terms of special values of 
$\zeta_W(s;{\frak g})$. This expression is called Witten's 
volume formula, which implies that 
\begin{align}
  \label{1-2}
    \zeta_W(2k;{\frak g})=C_{W}(2k,{\frak g})\pi^{2kn}
\end{align}
for any $k\in{\Bbb N}$, where 
$C_{W}(2k,{\frak g})\in{\Bbb Q}$ (see \cite[Theorem, p.506]{Za}). In general, the explicit value of $C_{W}(2k,{\frak g})$ was not determined in their work.

In our previous work \cite{KMT,KM3}, we defined the Bernoulli polynomials of root systems, and proved a formula which expresses $C_{W}(2k,{\frak g})$ in terms of those 
Bernoulli polynomials. Consequently we were able to obtain a certain generalization of \eqref{1-2}. We further gave some functional relations for zeta-functions of root systems which include \eqref{1-2} as special value-relations. In fact, we studied explicit functional relations for zeta-functions of $A_r$ type in \cite{MTF,KMT}, and of $B_r$ and $C_r$ types in \cite{KM3,KMTJC} (see also \cite{KMTpja}).

In the present paper we continue our research mentioned above. The main aim of the present paper is to study the zeta-function of $G_2$-type. 
In Section \ref{sec-2}, we define the Bernoulli polynomials of $G_2$-type and study the generating functions of them. By this consideration, we give \eqref{1-2} for $\zeta_2({\bf s};G_2)$ with explicit values of $C_W(2k,G_2)$ and more generalized results. In Section \ref{sec-3}, we consider analytic properties of $\zeta_2({\bf s};G_2)$ based on our previous paper \cite{KM2}. Actually we determine the possible singularities of $\zeta_2({\bf s};G_2)$. In Section \ref{sec-4}, we quote several lemmas which were shown in our previous papers \cite{KM3,KMTJC}. Furthermore we prove a certain analogue of them. These lemmas will play important roles in the next section. Finally, in Section \ref{sec-5}, by using these lemmas, we give explicit functional relations for $\zeta_2({\bf s};G_2)$, which include \eqref{1-2} at their special values. Recently, Zhao expressed several values $\zeta_2({\bf k};G_2)$ for ${\bf k}\in \mathbb{N}_0^6$ in terms of double polylogarithms and gave their approximate values in \cite{Zhao}. We will be able to give some of these values exactly (see Example \ref{Exam-2-1} and Remark \ref{Rem-4-7}). A part of these results has also been announced in our previous paper \cite{KMTpja}. 

Finally we remark that it is theoretically possible to prove the same type of results for zeta-functions of other exceptional types $E_6$, $E_7$, $E_8$ and $F_4$, by using our method. However, it may be considerably hard to apply our method actually to those cases, while it is interesting to determine $C_W(2k,{\frak g})$ explicitly in those cases. 

\bigskip

%%%%%%%%%%%%%%%%%%%%%%%%%%%%%%%%%%%%%%%%%%%%%%%%%%%%%%%%%%%%%%%
\section{Generating functions of the Bernoulli polynomials of $G_2$-type} \label{sec-2}
%%%%%%%%%%%%%%%%%%%%%%%%%%%%%%%%%%%%%%%%%%%%%%%%%%%%%%%%%%%%%%%

In our previous papers \cite{KM2,KM3}, we have already studied the general theory of zeta-functions of root systems. We apply it to the case of $G_2$ as follows.

% Let $V$ be a $2$-dimensional real vector space equipped with an inner
% product $\langle \cdot,\cdot\rangle$. 
% a $2$-dimensional real vector space equipped with an inner
% product $\langle \cdot,\cdot\rangle$. 
Let $\Delta=\Delta(G_2)\subset\mathfrak{h}^*$ be 
the root system %of
%a finite reduced root system %in $V$ 
of $G_2$-type and 
let $\mathfrak{h}_0=\mathbb{R}\alpha_1^\vee\oplus\mathbb{R}\alpha_2^\vee$
be a real vector subspace.
Let $\Delta_+$ and
$\Delta_-$ be the set of all positive roots and negative roots
respectively.  Then we have a decomposition of the root system
$\Delta=\Delta_+\coprod\Delta_-$.  We know that $\Delta_+$ is given by
\begin{equation}
  \Delta_+=\{\alpha_j\}_{j=1}^6, \label{2-1}
\end{equation}
where $\Psi=\{\alpha_1,\alpha_2\}$ is the set of fundamental roots and
\begin{equation}
\label{2-2}
  \begin{aligned}
    \alpha_3&=3\alpha_1+\alpha_2,\qquad &\alpha_3^\vee&=\alpha_1^\vee+\alpha_2^\vee\\
    \alpha_4&=3\alpha_1+2\alpha_2,\qquad &\alpha_4^\vee&=\alpha_1^\vee+2\alpha_2^\vee\\
    \alpha_5&=\alpha_1+\alpha_2,\qquad &\alpha_5^\vee&=\alpha_1^\vee+3\alpha_2^\vee\\
    \alpha_6&=2\alpha_1+\alpha_2,\qquad &\alpha_6^\vee&=2\alpha_1^\vee+3\alpha_2^\vee
  \end{aligned}
\end{equation}
(see, for example, Bourbaki \cite{Bourbaki}). 

\begin{figure}[h]
  \centering
\begin{overpic}[scale=0.8]{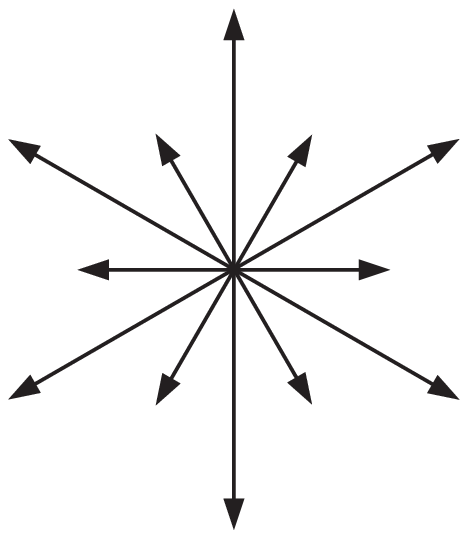}
  \put(73,48){$\alpha_1$}
  \put(0,76){$\alpha_2$}
  \put(82,76){$\alpha_3$}
  \put(40,98){$\alpha_4$}
  \put(26,76){$\alpha_5$}
  \put(56,76){$\alpha_6$}
\end{overpic}
\begin{overpic}[scale=0.8]{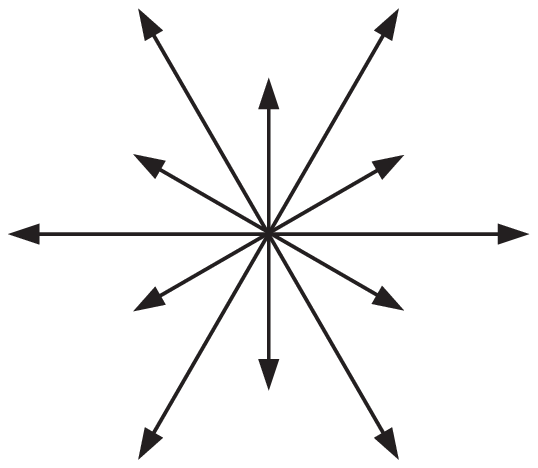}
  \put(96,43){$\alpha_1^\vee$}
  \put(20,62){$\alpha_2^\vee$}
  \put(72,62){$\alpha_3^\vee$}
  \put(48,74){$\alpha_4^\vee$}
  \put(24,86){$\alpha_5^\vee$}
  \put(70,86){$\alpha_6^\vee$}
\end{overpic}
  \caption{$G_2$}
  \label{fig:G2}
\end{figure}

Let $W=W(G_2)$ be the Weyl group of $G_2$-type.
For $w\in W$, we set 
\begin{equation}
\label{2-4}
\Delta_w=\Delta_+\cap w^{-1}\Delta_-.  
\end{equation}

From \eqref{1-3} and \eqref{2-2}, we see that the zeta-function of the root system of $G_2$-type can be given by 
\begin{equation}
\begin{split}
\zeta_2({\bf s};G_2)&  =\zeta_2(s_1,s_2,s_3,s_4,s_5,s_6;G_2)\\
& =  \sum_{m=1}^\infty\sum_{n=1}^\infty
  \frac{1}{m^{s_1}n^{s_2}(m+n)^{s_3}(m+2n)^{s_4}(m+3n)^{s_5}(2m+3n)^{s_6}},
\end{split}
\label{def-G2}
\end{equation}
where $s_j=s_{\alpha_j}$ ($1\leq j\leq 6$).
Furthermore, for $i=\sqrt{-1}$ and $\mathbf{y}=y_1\alpha_1^\vee+y_2\alpha_2^\vee\in \mathfrak{h}_0$, we define
\begin{equation}
\begin{split}
\zeta_2(\mathbf{s},\mathbf{y};G_2)&  =\zeta_2(s_1,s_2,s_3,s_4,s_5,s_6,y_1,y_2;G_2)\\
& =  \sum_{m=1}^\infty\sum_{n=1}^\infty
  \frac{e^{2\pi i (my_1+ny_2)}}{m^{s_1}n^{s_2}(m+n)^{s_3}(m+2n)^{s_4}(m+3n)^{s_5}(2m+3n)^{s_6}}.
\end{split}
\label{def-G2-y}
\end{equation}
Note that we only use the case when $\mathbf{y}=\mathbf{0}$ in this
paper. However, we study the case of general $\mathbf{y}$ here for the
convenience of our research in the future.  With the above notation
such as $W=W(G_2)$ and $\Delta=\Delta(G_2)$, we let
\begin{equation}
\label{eq:def-S}
  S(\mathbf{s},\mathbf{y};G_2)
=\sum_{w\in W}
    \Bigl(\prod_{\alpha\in\Delta_{w^{-1}}}(-1)^{-s_{\alpha}}\Bigr)
\zeta_2(w^{-1}\mathbf{s},w^{-1}\mathbf{y};G_2),
\end{equation}
where 
$(w^{-1}\mathbf{s})_\alpha=s_{w\alpha}$ with the identification $s_\alpha=s_{-\alpha}$
and $w^{-1}\mathbf{y}$ is the usual action of $\mathbf{y}$ by $w^{-1}$.

This $S(\mathbf{s},\mathbf{y};G_2)$
is a ``Weyl group symmetric'' linear combination of zeta-functions of root
systems, which plays a fundamental role in the 
study of value-relations and functional relations
in \cite{KM3}.

In order to evaluate $S(\mathbf{s},\mathbf{y};G_2)$ at positive
integers, we consider Bernoulli polynomials
$P(\mathbf{k},\mathbf{y};G_2)$ via their generating function $F({\bf
  t},{\bf y};G_2)$. This type of generalized Bernoulli polynomials
associated with any root system was first introduced in \cite{KM3},
and was further studied in \cite{KMT-L}.  For a real number $x$, let
$\{x\}$ denote its fractional part $x-[x]$.  Applying Theorem 4.1 in
\cite{KMT-L} to the case of $G_2$-type, we obtain
\begin{align*}
& F({\bf t},{\bf y};G_2)=F(t_1,t_2,t_3,t_4,t_5,t_6,y_1,y_2;G_2)\\
& =t_1t_2t_3t_4t_5t_6 \\
& \times\Bigg(
\frac{e^{\{y_1\} t_1+\{y_2\} t_2} }
{(e^{t_1}-1) (e^{t_2}-1) (t_1+t_2-t_3) (t_1+2 t_2-t_4) (t_1+3 t_2-t_5) (2 t_1+3 t_2-t_6)}
\\
&+\frac{e^{\{y_1-y_2\} t_1+\{y_2\} t_3} }
{(e^{t_1}-1) (e^{t_3}-1) (t_1+t_2-t_3) (t_1-2 t_3+t_4) (2 t_1-3 t_3+t_5) (t_1-3 t_3+t_6)}
\\
&+\frac{e^{\left\{y_1-\frac{y_2}{2}+\frac{1}{2}\right\} t_1+\left\{\frac{y_2}{2}+\frac{1}{2}\right\} t_4}
+e^{\left\{y_1-\frac{y_2}{2}\right\} t_1+\left\{\frac{y_2}{2}\right\} t_4}}
{2 (e^{t_1}-1) (e^{t_4}-1) \left(\frac{t_1}{2}+t_2-\frac{t_4}{2}\right) \left(\frac{t_1}{2}-t_3+\frac{t_4}{2}\right) \left(\frac{t_1}{2}-\frac{3 t_4}{2}+t_5\right) \left(\frac{t_1}{2}+\frac{3 t_4}{2}-t_6\right)}
\\
&-\frac{e^{\left\{y_1-\frac{y_2}{3}+\frac{2}{3}\right\} t_1+\left\{\frac{y_2}{3}+\frac{1}{3}\right\} t_5}
+e^{\left\{y_1-\frac{y_2}{3}+\frac{4}{3}\right\} t_1+\left\{\frac{y_2}{3}+\frac{2}{3}\right\} t_5}
+e^{\left\{y_1-\frac{y_2}{3}\right\} t_1+\left\{\frac{y_2}{3}\right\} t_5}}
{3 (e^{t_1}-1) (e^{t_5}-1) \left(\frac{t_1}{3}-t_4+\frac{2 t_5}{3}\right) \left(\frac{t_1}{3}+t_2-\frac{t_5}{3}\right) \left(\frac{2 t_1}{3}-t_3+\frac{t_5}{3}\right) (t_1+t_5-t_6)}
\\
&-\frac{e^{\left\{y_1-\frac{2 y_2}{3}+\frac{1}{3}\right\} t_1+\left\{\frac{y_2}{3}+\frac{1}{3}\right\} t_6}
+e^{\left\{y_1-\frac{2 y_2}{3}+\frac{2}{3}\right\} t_1+\left\{\frac{y_2}{3}+\frac{2}{3}\right\} t_6}
+e^{\left\{y_1-\frac{2 y_2}{3}\right\} t_1+\left\{\frac{y_2}{3}\right\} t_6}}
{3 (e^{t_1}-1) (e^{t_6}-1) (t_1+t_5-t_6) \left(\frac{t_1}{3}+t_4-\frac{2 t_6}{3}\right) \left(\frac{2 t_1}{3}+t_2-\frac{t_6}{3}\right) \left(\frac{t_1}{3}-t_3+\frac{t_6}{3}\right)}
\\
&-\frac{e^{(1-\{y_1-y_2\}) t_2+\{y_1\} t_3}}
{(e^{t_2}-1) (e^{t_3}-1) (t_1+t_2-t_3) (t_2+t_3-t_4) (2 t_2+t_3-t_5) (t_2+2 t_3-t_6)}
\\
&-\frac{e^{(1-\{2 y_1-y_2\}) t_2+\{y_1\} t_4} }
{(e^{t_2}-1) (e^{t_4}-1) (t_1+2 t_2-t_4) (t_2+t_3-t_4) (t_2+t_4-t_5) (t_2-2 t_4+t_6)}
\\
&+\frac{e^{(1-\{3 y_1-y_2\}) t_2+\{y_1\} t_5} }
{(e^{t_2}-1) (e^{t_5}-1) (t_1+3 t_2-t_5) (2 t_2+t_3-t_5) (t_2+t_4-t_5) (3 t_2-2 t_5+t_6)}
\\
&+\frac{e^{\left\{-\frac{3 y_1}{2}+y_2+\frac{1}{2}\right\} t_2+\left\{\frac{y_1}{2}+\frac{1}{2}\right\} t_6}+e^{\left(1-\left\{\frac{3 y_1}{2}-y_2\right\}\right) t_2+\left\{\frac{y_1}{2}\right\} t_6}}
{2 (e^{t_2}-1) (e^{t_6}-1) \left(t_1+\frac{3 t_2}{2}-\frac{t_6}{2}\right) \left(\frac{t_2}{2}+t_3-\frac{t_6}{2}\right) \left(\frac{t_2}{2}-t_4+\frac{t_6}{2}\right) \left(\frac{3 t_2}{2}-t_5+\frac{t_6}{2}\right)}
\\
&-\frac{e^{\{2 y_1-y_2\} t_3+(1-\{y_1-y_2\}) t_4}}
{(e^{t_3}-1) (e^{t_4}-1) (t_2+t_3-t_4) (t_1-2 t_3+t_4) (t_3-2 t_4+t_5) (t_3+t_4-t_6)}
\\
&+\frac{e^{\left\{\frac{3 y_1}{2}-\frac{y_2}{2}\right\} t_3+\left(1-\left\{\frac{y_1}{2}-\frac{y_2}{2}\right\}\right) t_5}}
{2 (e^{t_3}-1) (e^{t_5}-1) \left(t_2+\frac{t_3}{2}-\frac{t_5}{2}\right) \left(\frac{t_3}{2}-t_4+\frac{t_5}{2}\right) \left(t_1-\frac{3 t_3}{2}+\frac{t_5}{2}\right) \left(\frac{3 t_3}{2}+\frac{t_5}{2}-t_6\right)}
\\
&+\frac{e^{\{3 y_1-2 y_2\} t_3+(1-\{y_1-y_2\}) t_6}}
{(e^{t_3}-1) (e^{t_6}-1) (3 t_3+t_5-2 t_6) (t_2+2 t_3-t_6) (t_3+t_4-t_6) (t_1-3 t_3+t_6)}
\\
&-\frac{e^{\{3 y_1-y_2\} t_4+(1-\{2 y_1-y_2\}) t_5} }
{(e^{t_4}-1) (e^{t_5}-1) (t_2+t_4-t_5) (t_3-2 t_4+t_5) (t_1-3 t_4+2 t_5) (3 t_4-t_5-t_6)}
\\
&-\frac{e^{(1-\{3 y_1-2 y_2\}) t_4+\{2 y_1-y_2\} t_6} }
{(e^{t_4}-1) (e^{t_6}-1) (t_1+3 t_4-2 t_6) (t_3+t_4-t_6) (t_2-2 t_4+t_6) (3 t_4-t_5-t_6)}
\\
&+\frac{e^{\left(1-\left\{y_1-\frac{2 y_2}{3}\right\}\right) t_5+\left\{y_1-\frac{y_2}{3}\right\} t_6}
+e^{\left\{-y_1+\frac{2 y_2}{3}+\frac{2}{3}\right\} t_5+\left\{y_1-\frac{y_2}{3}+\frac{2}{3}\right\} t_6}
+e^{\left\{-y_1+\frac{2 y_2}{3}+\frac{4}{3}\right\} t_5+\left\{y_1-\frac{y_2}{3}+\frac{4}{3}\right\} t_6}}
{3 (e^{t_5}-1) (e^{t_6}-1) (t_1+t_5-t_6) \left(t_3+\frac{t_5}{3}-\frac{2 t_6}{3}\right) \left(t_4-\frac{t_5}{3}-\frac{t_6}{3}\right) \left(t_2-\frac{2
   t_5}{3}+\frac{t_6}{3}\right)}\Bigg).
\end{align*}
Then $F({\bf t},{\bf y};G_2)$ is holomorphic at the origin and can be expanded as 
\begin{align*}
\label{eq:def_F}
  F(\mathbf{t},\mathbf{y};G_2)&=
  \sum_{\mathbf{k}\in \mathbb{N}_0^{6}}P(\mathbf{k},\mathbf{y};G_2)
  \prod_{\alpha\in\Delta_+}
  \frac{t_\alpha^{k_\alpha}}{k_\alpha!}
\end{align*}
for $\mathbf{y}\in \mathfrak{h}_0$. From our previous results \cite[Theorem 4.4,\,(4.19) and (4.20)]{KM3}, we obtain the following.

\begin{theorem} \label{T-2-1}
For $\mathbf{k}\in\mathbb{N}_0^{6}$,
\begin{equation}
\label{eq:S_B}
S(\mathbf{k},\mathbf{y};G_2)=
  \biggl(\prod_{\alpha\in\Delta_+}
  \frac{(2\pi i)^{k_\alpha}}{k_\alpha!}\biggr)
P(\mathbf{k},\mathbf{y};G_2).
\end{equation}
\end{theorem}

\begin{example}\label{Exam-2-1}
In the case $\mathbf{k}=(2m,2m,\ldots,2m)$ for $m\in \mathbb{N}$ and ${\bf y}={\bf 0}$, we see that 
\begin{equation}
S((2m),\mathbf{0};G_2)=12\zeta_2((2m);G_2). \label{S-value}
\end{equation}
On the other hand, from the definition of $F(\mathbf{t},\mathbf{y};G_2)$, we 
can calculate $P(\mathbf{k},\mathbf{0};G_2)$. 
Combining this fact with \eqref{S-value}, we can obtain the explicit values of $\zeta_2((2m);G_2)$, for example, 
\begin{align*}
& \zeta_2(2,2,2,2,2,2;G_2) = \frac{23}{297904566960} \pi^{12}; \\
& \zeta_2(4,4,4,4,4,4;G_2) = \frac{8165653}{1445838676129559305994400000} \pi^{24} ;\\
& \zeta_2(6,6,6,6,6,6;G_2) = \frac{55940539974690617}{131888156302530666544150214880458495963616000000} \pi^{36};\\
& \zeta_2(8,8,8,8,8,8;G_2) \\
& = \frac{47346365461279256768015189}{1485697621623958244738368714652675148113575302412190275200000000000} \pi^{48}.
\end{align*}
Furthermore, in the case when ${\bf k}=(2p,2q,2q,2q,2p,2p)$ $(p,q \in \mathbb{N})$, we have
$$S(\mathbf{k},\mathbf{0};\Delta)=12\zeta_2(2p,2q,2q,2q,2p,2p;G_2).$$
This is because the lengths of $\alpha_1$, $\alpha_5$ and $\alpha_6$ 
(and of $\alpha_2$, $\alpha_3$ and $\alpha_4$) are the same,
and the roots of the same length form a single Weyl-group orbit.
Hence we can obtain, for example,
\begin{align*}
& \zeta_2(2,4,4,4,2,2;G_2) = \frac{467}{213955059990672000} \pi^{18},\\ 
& \zeta_2(4,2,2,2,4,4;G_2) 
  = \frac{20771}{106061802338575923840} \pi^{18},  \\
& \zeta_2(2,6,6,6,2,2;G_2) 
  = \frac{91027}{1449347623006311428400000} \pi^{24},\\
& \zeta_2(6,2,2,2,6,6;G_2) 
  = \frac{391420483}{770242750118097151820324400000} \pi^{24},\\
& \zeta_2(2,8,8,8,2,2;G_2)  = \frac{19152444887}{10564558460425628849656425960000000} \pi^{30},\\
& \zeta_2(8,2,2,2,8,8;G_2)  = \frac{1802533972626341}{1364308801602394759022133342471831480000000} \pi^{30}.
\end{align*}
It is possible to compute the numerical values of the left-hand sides of the above formulas from the definition \eqref{def-G2}.
We have already checked that those numerical values agree with the above formulas.
\end{example}

\bigskip

%%%%%%%%%%%%%%%%%%%%%%%%%%%%%%%%%%%%%%%%%%%%%%%%%%%%%%%%%%%%%%%
\section{Analytic properties of the zeta-function of $G_2$-type} \label{sec-3}
%%%%%%%%%%%%%%%%%%%%%%%%%%%%%%%%%%%%%%%%%%%%%%%%%%%%%%%%%%%%%%%

In the preceding section we studied ``value-relations'' for $\zeta_2({\bf s};G_2)$, but we can further discuss ``functional relations'' for this function. For this purpose, we first consider analytic properties. 

\begin{theorem} \label{T-3-1}
The function $\zeta_2(s_1,s_2,s_3,s_4,s_5,s_6;G_2)$ can be continued meromorphically to the whole space $\mathbb{C}^6$, and its possible singularities are located on the subsets of $\mathbb{C}^6$ defined by one of the equations:
\begin{align*}
& s_1+s_3+s_4+s_5+s_6=1-l\quad (l\in \mathbb{N}_0),\\
& s_2+s_3+s_4+s_5+s_6=1-l\quad (l\in \mathbb{N}_0),\\
& s_1+s_2+s_3+s_4+s_5+s_6=2.
\end{align*}
\end{theorem}

\smallskip

The meromorphic continuation of $\zeta_2({\bf s};G_2)$, as well as zeta-functions of other exceptional algebras, can be deduced from earlier results given by Essouabri \cite{Ess95,Ess97}, the second-named author \cite[Theorem 3]{Ma2}, and the first-named author \cite{Komori09}. However the following argument of determining the possible singularities also includes a proof of meromorphic continuation. 

At the end of \cite[Remark 6.4]{KM2}, we discussed when the determination of possible singularities can be achieved just by shifting the path of integration. The arrow from $G_2$ to $C_2$ in the diagram in \cite[Section 5]{KM2} is horizontal, hence this is the case when the shifting of the path is sufficient. Therefore our argument here is not so complicated, similar to that in \cite{Ma,Ma0,Ma1}. Note that such simple shifting argument is not sufficient when one studies analytic properties of  zeta-functions of other exceptional algebras.

At first, assume $\Re s_j$ $(1\leq j\leq 6)$ are sufficiently large. The Mellin-Barnes integral formula is 
\begin{equation}
(1+\lambda)^{-s}=\frac{1}{2\pi i}\int_{(c)}\frac{\Gamma(s+z)\Gamma(-z)}{\Gamma(s)}\lambda^z dz, \label{3-1}
\end{equation}
where $\Re s>0$, $|\textrm{arg}\lambda|<\pi$, $\lambda \not=0$, $c \in \mathbb{R}$ with $-\Re s<c<0$ and the path $(c)$ of integration is the vertical line $\Re z=c$. By using \eqref{3-1}, we first prove an integral expression of $\zeta_2({\bf s};G_2)$ in terms of the zeta-function of $C_2$-type defined by 
$$\zeta_2(s_1,s_2,s_3,s_4;C_2)=\sum_{m=1}^\infty \sum_{n=1}^\infty \frac{1}{m^{s_1}n^{s_2}(m+n)^{s_3}(m+2n)^{s_4}}$$
(see \cite[(6.1)]{KM2} and \cite[Example 7.3]{KM3}). Writing
$$(m+3n)^{-s_5}=(m+2n)^{-s_5}\left(1+\frac{n}{m+2n}\right)^{-s_5}$$
and applying \eqref{3-1} to the second factor of the right-hand side, we have
\begin{equation*}
(m+3n)^{-s_5}=\frac{1}{2\pi i}\int_{(c_1)}\frac{\Gamma(s_5+z_1)\Gamma(-z_1)}{\Gamma(s_5)}(m+2n)^{-s_5-z_1}n^{z_1} dz_1, 
\end{equation*}
with $-\Re s_5<c_1<0$. Similarly 
\begin{equation*}
(2m+3n)^{-s_6}=\frac{1}{2\pi i}\int_{(c_2)}\frac{\Gamma(s_6+z_2)\Gamma(-z_2)}{\Gamma(s_6)}(m+2n)^{-s_6-z_2}(m+n)^{z_2} dz_2, 
\end{equation*}
with $-\Re s_6<c_2<0$. Hence
\begin{equation}
\begin{split}
\zeta_2({\bf s};G_2) & = \frac{1}{(2\pi i)^2}\int_{(c_1)}\int_{(c_2)}\frac{\Gamma(s_5+z_1)\Gamma(s_6+z_2)\Gamma(-z_1)\Gamma(-z_2)}{\Gamma(s_5)\Gamma(s_6)} \\
& \quad \times \zeta_2(s_1,s_2-z_1,s_3-z_2,s_4+s_5+s_6+z_1+z_2;C_2)dz_2 dz_1. 
\end{split}
\label{3-2}
\end{equation}
The singularities of $\zeta_2({\bf s};C_2)$ are determined by \cite[Theorem 6.2]{KM2}. We find that the singularities of the zeta factor on the right-hand side of \eqref{3-2} is 
\begin{align}
& s_1+s_3+s_4+s_5+s_6+z_1=1-l\quad (l\in \mathbb{N}_0),  \label{3-3}\\
& s_2+s_3+s_4+s_5+s_6=1-l\quad (l\in \mathbb{N}_0), \label{3-4}\\
& s_1+s_2+s_3+s_4+s_5+s_6=2. \label{3-5}
\end{align}
Let $L$ be a large positive integer, and define 
\begin{align*}
& \Phi({\bf s}) = \prod_{l=0}^{L-1}(s_2+s_3+s_4+s_5+s_6-1+l)(s_1+s_2+s_3+s_4+s_5+s_6-2).
\end{align*}
We can rewrite \eqref{3-2} as 
\begin{equation}
\begin{split}
\zeta_2({\bf s};G_2) & = \Phi({\bf s})^{-1}\frac{1}{2\pi i}\int_{(c_2)}\frac{\Gamma(s_6+z_2)\Gamma(-z_2)}{\Gamma(s_6)}I({\bf s},z_2)dz_2, 
\end{split}
\label{3-6}
\end{equation}
where 
\begin{equation}
\begin{split}
I({\bf s},z_2) & =\frac{1}{2\pi i}\int_{(c_1)}\frac{\Gamma(s_5+z_1)\Gamma(-z_1)}{\Gamma(s_5)}\Phi({\bf s})\\
& \quad \times \zeta_2(s_1,s_2-z_1,s_3-z_2,s_4+s_5+s_6+z_1+z_2;C_2)dz_1.
\end{split}
\label{3-7}
\end{equation}
The singularities of the integrand on the right-hand side of \eqref{3-7} are $z_1=-s_5-l$, $z_1=l$ $(l\in \mathbb{N}_0)$, singularities of type \eqref{3-3}, and of type \eqref{3-4} with $l\geq L$. 
Shifting the path to $\Re z_1=M_1-\varepsilon$, where $M_1$ is a large positive integer and $\varepsilon$ is a small positive number, and counting the residues at $z_1=0,1,2,\ldots,M_1-1$, we obtain
\begin{equation}
\begin{split}
I({\bf s},z_2) & =\sum_{m_1=0}^{M_1-1}\binom{-s_5}{m_1}\Phi({\bf s})\zeta_2(s_1,s_2-m_1,s_3-z_2,s_4+s_5+s_6+m_1+z_2;C_2)\\
& +\frac{1}{2\pi i}\int_{(M_1-\varepsilon)}\frac{\Gamma(s_5+z_1)\Gamma(-z_1)}{\Gamma(s_5)}\Phi({\bf s})\\
& \quad \times \zeta_2(s_1,s_2-z_1,s_3-z_2,s_4+s_5+s_6+z_1+z_2;C_2)dz_1,
\end{split}
\label{3-8}
\end{equation}
where 
\begin{equation*}
\binom{s}{k}=
\begin{cases} 
\frac{s(s-1)(s-2)\cdots (s-k+1)}{k!} & (k\in \mathbb{N}),\\
\qquad 1 & (k=0).
\end{cases}
\end{equation*}
The above integral is holomorphic in the region
\begin{align*}
& \Re s_5>-M_1+\varepsilon, \\
& \Re (s_1+s_3+s_4+s_5+s_6)>-M_1+1+\varepsilon,
\end{align*}
and 
\begin{equation}
\Re (s_2+s_3+s_4+s_5+s_6)>1-L. \label{3-9}
\end{equation}
Since $M_1$ is arbitrary, we now find that $I({\bf s},z_2)$ is continued meromorphically to the region \eqref{3-9}. 
We can also show that $I({\bf s},z_2)$ is of polynomial order with respect to the imaginary parts of variables. The singularities of $I({\bf s},z_2)$ in the region \eqref{3-9} are located only on 
\begin{equation}
s_1+s_3+s_4+s_5+s_6=1-l\quad (l\in \mathbb{N}_0), \label{3-10}
\end{equation}
which are coming from the zeta-factors in the sum-part on the right-hand side of \eqref{3-8}. 

Now go back to the situation when $\Re s_j$ $(1\leq j \leq 6)$ are large, and consider \eqref{3-6}. Let 
\begin{align*}
& \Psi({\bf s}) = \prod_{l=0}^{L-1}(s_1+s_3+s_4+s_5+s_6-1+l),
\end{align*}
insert $\Phi({\bf s})^{-1}\Psi({\bf s})^{-1}$ on the right-hand side of \eqref{3-6}, and shift the path to $\Re z_2=M_2-\varepsilon$ to obtain 
\begin{equation}
\begin{split}
\zeta_2({\bf s};G_2) & = \Phi({\bf s})^{-1}\Psi({\bf s})^{-1}\bigg\{ \sum_{m_2=0}^{M_2-1}\binom{-s_6}{m_2}\Psi({\bf s})I({\bf s},m_2) \\
& \quad\quad\quad +\frac{1}{2\pi i}\int_{(M_2-\varepsilon)}\frac{\Gamma(s_6+z_2)\Gamma(-z_2)}{\Gamma(s_6)}\Psi({\bf s})I({\bf s},z_2)dz_2\bigg\}. 
\end{split}
\label{3-11}
\end{equation}
This gives the continuation of $\zeta_2({\bf s};G_2)$ to the region \eqref{3-9}. Since $L$ is arbitrary, we obtain the meromorphic continuation of $\zeta_2({\bf s};G_2)$ to the whole space. Its possible singularities are coming from $\Phi({\bf s})^{-1}\Psi({\bf s})^{-1}$ and $\Psi({\bf s})I({\bf s},m_2)$ on the right-hand side of \eqref{3-11}, which are exactly those stated in the theorem. 
This completes the proof of Theorem \ref{T-3-1}.

\bigskip

%%%%%%%%%%%%%%%%%%%%%%%%%%%%%%%%%%%%%%%%%%%%%%%%%%%%%%%%%%%%%%%
\section{Preliminary Lemmas} \label{sec-4}
%%%%%%%%%%%%%%%%%%%%%%%%%%%%%%%%%%%%%%%%%%%%%%%%%%%%%%%%%%%%%%%

In this section, we quote several lemmas from our previous papers \cite{MNO,TsC,KMTJC} and further prove an analogue of them. These will play important roles in the next section. From now on, the symbol $\{\ \}$ implies ordinary curly parentheses, not the fractional part. 

\bigskip

\begin{lemma}[\cite{MNO}\ Lemma 2.1]  \label{L-4-1} 
Let $\phi(s):=\sum_{n \geq 1}(-1)^n n^{-s}=\left( 2^{1-s}-1 \right)\zeta(s)$, and $f,g\,:\,\mathbb{N}_0 \to \mathbb{C}$ be arbitrary functions. Then, for $a \in \mathbb{N}$, 
\begin{align}
& \sum_{k=0}^{a} \phi(a-k) \lambda_{a-k} \sum_{\mu=0}^{[k/2]} f(k-2\mu)\frac{(i\pi)^{2\mu}}{(2\mu)!} =\sum_{\xi=0}^{[a/2]}\zeta(2\xi)f(a-2\xi), \label{eq-4-1}
\end{align}
and
\begin{align}
& \ \sum_{k=1}^{a} \phi(a-k) \lambda_{a-k} \sum_{\mu=0}^{[(k-1)/2]} g(k-2\mu)\frac{(i\pi)^{2\mu}}{(2\mu+1)!} =-\frac{1}{2}g(a), \label{eq-4-2}
\end{align}
where $\lambda_\nu:=(1+(-1)^\nu)/2$ for $\nu \in \mathbb{Z}$. 
\end{lemma}

\begin{lemma}[\cite{TsC}\ Lemma 4.4] \label{L-4-2}
Let $\{ P_{2h}\},\ \{Q_{2h}\},\ \{R_{2h}\}$ be sequences such that
\begin{equation*}
P_{2h}=\sum_{j=0}^{h}R_{2h-2j} \frac{(i\pi)^{2j}}{(2j)!},\ Q_{2h}=\sum_{j=0}^{h}R_{2h-2j}\frac{(i\pi)^{2j}}{(2j+1)!}
\end{equation*}
for any $h \in \mathbb{N}_0$. Then  
\begin{align}
& P_{2h}=-2\sum_{\tau=0}^{h}\zeta(2h-2\tau)Q_{2\tau}, \label{P-2h} \\
& Q_{2h}=\frac{2}{\pi^2}\sum_{\tau=0}^{h}\left(2^{2h-2\tau+2}-1\right)\zeta(2h-2\tau+2)P_{2\tau} \label{Q-2h}
\end{align}
for any $h \in \mathbb{N}_0$. 
\end{lemma}

Note that, in \cite[Lemma 4.4]{TsC}, we proved only \eqref{P-2h}. However, by just the same method, we can easily obtain \eqref{Q-2h}.

\begin{lemma}[\cite{KMTJC}\ Lemma 6.3] \label{L-4-3} Let $h \in \mathbb{N}$, and
\begin{align*}
& {\frak C}:=\left\{ C(l) \in \mathbb{C}\,|\, l \in \mathbb{Z},\ l \not=0 \right\}, \\
& {\frak D}:=\left\{ D(N;m;\eta) \in \mathbb{R}\,|\, N,m,\eta \in \mathbb{Z},\ N \not=0,\ m \geq 0,\ 1 \leq \eta \leq h\right\}, \\
& {\frak A}:=\{ a_\eta \in \mathbb{N}\,|\,1 \leq \eta \leq h\}
\end{align*}
be sets of numbers indexed by integers. Assume that the infinite series appearing in 
\begin{align}
\sum_{N \in \mathbb{Z} \atop N \not=0}(-1)^{N}C(N)e^{iN\theta} & -2\sum_{\eta=1}^{h}\sum_{k=0}^{a_\eta}\phi(a_\eta-k)\lambda_{a_\eta-k} \label{eq-4-3} \\
& \ \ \times \sum_{\xi=0}^{k}\left\{ \sum_{N \in \mathbb{Z} \atop N \not=0}(-1)^N D(N;k-\xi;\eta)e^{iN\theta}\right\}\frac{(i\theta)^\xi}{\xi!} \notag
\end{align}
are absolutely convergent for $\theta \in [-\pi,\pi]$, and that (\ref{eq-4-3}) is a constant function for $\theta \in [-\pi,\pi]$. Then, for $d \in \mathbb{N}_0$, 
\begin{align}
& \sum_{N \in \mathbb{Z} \atop N \not= 0}\frac{(-1)^{N}C(N)e^{iN\theta}}{N^d} -2\sum_{\eta=1}^{h}\sum_{k=0}^{a_\eta}\phi(a_\eta-k)\lambda_{a_\eta-k} \label{eq-4-4} \\
& \ \ \ \ \ \times \sum_{\xi=0}^{k}\bigg\{ \sum_{\omega=0}^{k-\xi}\binom{\omega+d-1}{\omega}(-1)^{\omega} \sum_{m \in \mathbb{Z} \atop m \not= 0}\frac{(-1)^m D(m;k-\xi-\omega;\eta)e^{im\theta}}{m^{d+\omega}}\bigg\}\frac{(i\theta)^\xi}{\xi!} \notag \\
& \ +2\sum_{k=0}^{d}\phi(d-k)\lambda_{d-k} \sum_{\xi=0}^{k}\bigg\{ \sum_{\eta=1}^{h} \sum_{\omega=0}^{a_\eta-1}\binom{\omega+k-\xi}{\omega}(-1)^{\omega}\notag \\
& \hspace{1in} \times \sum_{m \in \mathbb{Z} \atop m \not=0}\frac{D(m;a_\eta-1-\omega;\eta)}{m^{k-\xi+\omega+1}}\bigg\}\frac{(i\theta)^\xi}{\xi!} =0 \notag
\end{align}
holds for $\theta \in [-\pi,\pi]$, where the infinite series appearing on the left-hand side of (\ref{eq-4-4}) are absolutely convergent for $\theta \in [-\pi,\pi]$.
\end{lemma}

Now we prepare the following lemma which is an analogue of Lemma \ref{L-4-3}. 

\begin{lemma} \label{L-4-4} Let $h \in \mathbb{N}$, 
\begin{align*}
& {\mathcal{A}}:=\left\{ \alpha(l) \in \mathbb{C}\,|\, l \in \mathbb{Z},\ l \not=0 \right\}, \\
& {\mathcal{B}}:=\left\{ \beta(N;m;\eta) \in \mathbb{R}\,|\, N,m,\eta \in \mathbb{Z},\ N \not=0,\ m \geq 0,\ 1 \leq \eta \leq h\right\}, \\
& {\mathcal{C}}:=\{ c_\eta \in \mathbb{N}\,|\,1 \leq \eta \leq h\}
\end{align*}
be sets of numbers indexed by integers, and 
\begin{equation}
\begin{split}
S_{\pm}(\theta)=\sum_{m \in \mathbb{Z} \atop m \not=0}(\pm i)^{m}\alpha(m)e^{im\theta/2} & -2\sum_{\eta=1}^{h}\sum_{k=0}^{c_\eta}\phi(c_\eta-k)\lambda_{c_\eta-k}  \\
& \ \ \times \sum_{\xi=0}^{k}\left\{ \sum_{m \in \mathbb{Z} \atop m \not=0}(\pm i)^m \beta(m;k-\xi;\eta)e^{im\theta/2}\right\}\frac{(i\theta)^\xi}{\xi!}.
\end{split}
\label{eq-4-10}
\end{equation}
Assume that both of the right-hand sides of $S_{\pm}(\theta)$ in \eqref{eq-4-10} are absolutely convergent for $\theta \in [-\pi,\pi]$, and that both $S_{+}(\theta)$ and $S_{-}(\theta)$ are constant functions on $[-\pi,\pi]$. Then, for $d \in \mathbb{N}$, 
\begin{equation}
\begin{split}
& \sum_{m \in \mathbb{Z} \atop m \not= 0}\frac{\alpha(m)}{m^{2d}} -2\sum_{\eta=1}^{h}\sum_{k=0}^{[c_\eta/2]}\zeta(2k) \sum_{\omega=0}^{c_\eta-2k}\binom{\omega+2d-1}{\omega}(-2)^{\omega} \sum_{m \in \mathbb{Z} \atop m \not= 0}\frac{\beta(m;c_\eta-2k-\omega;\eta)}{m^{2d+\omega}} \\
& \ +2\sum_{k=0}^{d}\zeta(2k)2^{-2k} \sum_{\eta=1}^{h} \sum_{\omega=0}^{c_\eta-1}\binom{\omega+2d-2k}{\omega}(-2)^{\omega} \\
& \hspace{1in} \times \sum_{m \in \mathbb{Z} \atop m \not=0}\frac{((-1)^m+1)\beta(m;c_\eta-1-\omega;\eta)}{m^{2d-2k+\omega+1}} \\
& \ -2\sum_{k=0}^{d}\zeta(2k)\left(1-2^{-2k}\right) \sum_{\eta=1}^{h} \sum_{\omega=0}^{c_\eta-1}\binom{\omega+2d-2k}{\omega}(-2)^{\omega} \\
& \hspace{1in} \times \sum_{m \in \mathbb{Z} \atop m \not=0}\frac{((-1)^m-1)\beta(m;c_\eta-1-\omega;\eta)}{m^{2d-2k+\omega+1}}=0 
\end{split}
\label{eq-4-11} 
\end{equation}
for $\theta \in [-\pi,\pi]$, where the infinite series appearing on the left-hand side of (\ref{eq-4-11}) are absolutely convergent for $\theta \in [-\pi,\pi]$.\end{lemma}

\begin{proof}
Put
\begin{equation}
\begin{split}
\mathcal{G}_N^\pm(\theta)& =\mathcal{G}_N^\pm(\theta;\mathcal{A};\mathcal{B};\mathcal{C}) \\
& :=\left( \frac{2}{i}\right)^N \bigg\{ \sum_{m \in \mathbb{Z} \atop m \not= 0}\frac{\left(i^m \pm i^{-m}\right)\alpha(m)e^{im\theta/2}}{m^N} \\
& \ \ -2\sum_{\eta=1}^{h}\sum_{j=0}^{c_\eta}\phi(c_\eta-j)\lambda_{c_\eta-j} \sum_{\rho=0}^{j} \sum_{\omega=0}^{j-\rho} \binom{N-1+\omega}{\omega}(-2)^{\omega}\\
& \ \ \times \sum_{m \in \mathbb{Z} \atop m \not= 0}\frac{\left(i^m \pm i^{-m}\right)\beta(m;j-\rho-\omega;\eta)e^{im\theta/2}}{m^{N+\omega}}\frac{(i\theta)^\rho}{\rho!}\bigg\}\quad \ (N\in \mathbb{N}_0).
\end{split}
\label{eq-4-12}
\end{equation}
From the assumption, we see that $\mathcal{G}_0^\pm(\theta)$ are constant functions for $\theta \in [-\pi,\pi]$. Also, by 
using the relation
$$-\binom{m-1}{l-1}+\binom{m}{l}=\binom{m-1}{l}\ \ \ (l,m \in \mathbb{N}),$$
we can check that 
\begin{equation}
\frac{d}{d\theta}\mathcal{G}_N^\pm(\theta)=\mathcal{G}_{N-1}^\pm(\theta)\ \ (N \in \mathbb{N}). \label{eq-4-13}
\end{equation}
Repeating the indefinite integration, we can write
\begin{equation}
\left( \frac{i}{2}\right)^N \mathcal{G}_N^\pm(\theta)=\sum_{j=0}^{N}\mathfrak{C}_{N-j}^\pm \frac{(i\theta/2)^j}{j!}\ \ (N \in \mathbb{N}_0) \label{eq-4-14}
\end{equation}
for some $\{\mathfrak{C}_n^\pm \in \mathbb{C}\,|\,n \in \mathbb{N}_0\}$. 
Putting $N=2d+1$ for $d \in \mathbb{N}$ and $\theta=\pi$ in \eqref{eq-4-14}, we obtain 
\begin{equation}
\frac{(-1)^d}{2\pi}\left\{ \mathcal{G}_{2d+1}^+(\pi)-\mathcal{G}_{2d+1}^+(-\pi)\right\}=\sum_{\nu=0}^{d}\mathfrak{C}_{2d-2\nu}^+ 2^{2d-2\nu}\frac{(i\pi)^{2\nu}}{(2\nu+1)!}. \label{eq-4-15}
\end{equation}
Similarly, putting $N=2d$ and $\theta=\pi$ in \eqref{eq-4-14}, we have
\begin{equation}
\frac{(-1)^d}{2}\left\{ \mathcal{G}_{2d}^+(\pi)+\mathcal{G}_{2d}^+(-\pi)\right\}=\sum_{\nu=0}^{d}\mathfrak{C}_{2d-2\nu}^+ 2^{2d-2\nu}\frac{(i\pi)^{2\nu}}{(2\nu)!}. \label{eq-4-16}
\end{equation}
By Lemma \ref{L-4-2}, we have
\begin{equation}
\begin{split}
& \frac{(-1)^d}{2}\left\{ \mathcal{G}_{2d}^+(\pi)+\mathcal{G}_{2d}^+(-\pi)\right\} \\
& \ \ =-\frac{1}{\pi}\sum_{\tau=0}^{d}\zeta(2d-2\tau)(-1)^\tau \left\{ \mathcal{G}_{2\tau+1}^+(\pi)-\mathcal{G}_{2\tau+1}^+(-\pi)\right\}. 
\end{split}
\label{eq-4-17}
\end{equation}
We will calculate each side of \eqref{eq-4-17} explicitly as follows. Note that $$\left( i^m\pm i^{-m}\right)^2=2\left\{ (-1)^m\pm 1\right\},\ \ \left( i^m+i^{-m}\right)\left( i^m-i^{-m}\right)=0.$$
By using \eqref{eq-4-1}, we have
\begin{equation}
\begin{split}
\mathcal{G}_{2d}^+(\pi) + \mathcal{G}_{2d}^+(-\pi) & =(-1)^d 2^{2d} \bigg\{ 2\sum_{m \in \mathbb{Z} \atop m \not= 0}\frac{\left((-1)^m+ 1\right)\alpha(m)}{m^{2d}} \\
& \ \ -4\sum_{\eta=1}^{h}\sum_{j=0}^{c_\eta}\phi(c_\eta-j)\lambda_{c_\eta-j} \sum_{\mu=0}^{[j/2]} \sum_{\omega=0}^{j-2\mu} \binom{2d-1+\omega}{\omega}(-2)^{\omega}\\
& \ \ \ \ \times \sum_{m \in \mathbb{Z} \atop m \not= 0}\frac{\left((-1)^m +1\right)\beta(m;j-2\mu-\omega;\eta)}{m^{2d+\omega}}\frac{(i\pi)^{2\mu}}{(2\mu)!}\bigg\}\\
& =(-1)^d 2^{2d+1} \bigg\{ \sum_{m \in \mathbb{Z} \atop m \not= 0}\frac{\left((-1)^m+ 1\right)\alpha(m)}{m^{2d}} \\
& \ \ -2\sum_{\eta=1}^{h}\sum_{\xi=0}^{[c_\eta/2]}\zeta(2\xi) \sum_{\omega=0}^{c_\eta-2\xi} \binom{2d-1+\omega}{\omega}(-2)^{\omega}\\
& \ \ \ \ \times \sum_{m \in \mathbb{Z} \atop m \not= 0}\frac{\left((-1)^m +1\right)\beta(m;c_\eta-2\xi-w;\eta)}{m^{2d+\omega}}\bigg\}.
\end{split}
\label{eq-4-18}
\end{equation}
Similarly, by using \eqref{eq-4-2}, we have
\begin{equation}
\begin{split}
& \mathcal{G}_{2d+1}^+(\pi) - \mathcal{G}_{2d+1}^+(-\pi) \\
& =(-1)^d 2^{2d+2}\pi \sum_{\eta=1}^{h}\sum_{\omega=0}^{c_\eta-1}\binom{2d+\omega}{\omega}(-2)^{\omega}\sum_{m \in \mathbb{Z} \atop m \not= 0}\frac{\left((-1)^m +1\right)\beta(m;c_\eta-1-w;\eta)}{m^{2d+1+\omega}}.
\end{split}
\label{eq-4-19}
\end{equation}
Substituting \eqref{eq-4-18} and \eqref{eq-4-19} into \eqref{eq-4-17}, we have
\begin{equation}
\begin{split}
\sum_{m \in \mathbb{Z} \atop m \not= 0}\frac{\left((-1)^m+ 1\right)\alpha(m)}{m^{2d}} & -2\sum_{\eta=1}^{h}\sum_{\xi=0}^{[c_\eta/2]}\zeta(2\xi) \sum_{\omega=0}^{c_\eta-2\xi} \binom{2d-1+\omega}{\omega}(-2)^{\omega}\\
& \ \ \ \ \times \sum_{m \in \mathbb{Z} \atop m \not= 0}\frac{\left((-1)^m +1\right)\beta(m;c_\eta-2\xi-w;\eta)}{m^{2d+\omega}}\\
& \ \ +4\sum_{\eta=1}^{h}\sum_{\xi=0}^{d}\sum_{\omega=0}^{c_\eta-1}\zeta(2\xi)2^{-2\xi}\binom{2d-2\xi+\omega}{\omega}(-2)^{\omega} \\
& \ \ \ \ \times \sum_{m \in \mathbb{Z} \atop m \not= 0}\frac{\left((-1)^m +1\right)\beta(m;c_\eta-1-w;\eta)}{m^{2d-2\xi+1+\omega}}=0.
\end{split}
\label{eq-4-20}
\end{equation}
Similarly to \eqref{eq-4-17}, we obtain
\begin{equation}
\begin{split}
& \frac{(-1)^d}{2}\left\{ \mathcal{G}_{2d}^-(\pi)-\mathcal{G}_{2d}^-(-\pi)\right\} \\
& \ \ =-\frac{1}{\pi}\sum_{\tau=0}^{d}\left(2^{2d-2\tau}-1\right)\zeta(2d-2\tau)(-1)^\tau \left\{ \mathcal{G}_{2\tau+1}^-(\pi)+\mathcal{G}_{2\tau+1}^-(-\pi)\right\}. 
\end{split}
\label{eq-4-21}
\end{equation}
Hence, as well as \eqref{eq-4-20}, we have
\begin{equation}
\begin{split}
\sum_{m \in \mathbb{Z} \atop m \not= 0}\frac{\left((-1)^m- 1\right)\alpha(m)}{m^{2d}} & -2\sum_{\eta=1}^{h}\sum_{\xi=0}^{[c_\eta/2]}\zeta(2\xi) \sum_{\omega=0}^{c_\eta-2\xi} \binom{2d-1+\omega}{\omega}(-2)^{\omega}\\
& \ \ \ \ \times \sum_{m \in \mathbb{Z} \atop m \not= 0}\frac{\left((-1)^m -1\right)\beta(m;c_\eta-2\xi-w;\eta)}{m^{2d+\omega}}\\
& \ \ +4\sum_{\eta=1}^{h}\sum_{\xi=0}^{d}\sum_{\omega=0}^{c_\eta-1}\zeta(2\xi)\left(1-2^{-2\xi}\right)\binom{2d-2\xi+\omega}{\omega}(-2)^{\omega} \\
& \ \ \ \ \times \sum_{m \in \mathbb{Z} \atop m \not= 0}\frac{\left((-1)^m -1\right)\beta(m;c_\eta-1-w;\eta)}{m^{2d-2\xi+1+\omega}}=0.
\end{split}
\label{eq-4-22}
\end{equation}
Combining \eqref{eq-4-20} and \eqref{eq-4-22}, we obtain \eqref{eq-4-11}.
\end{proof}

\ 

%%%%%%%%%%%%%%%%%%%%%%%%%%%%%%%%%%%%%%%%%%%%%%%%%%%%%%%%%%%%%%%%%%
\section{Functional relations for $\zeta_2({\boldsymbol{s}};G_2)$}
\label{sec-5}
%%%%%%%%%%%%%%%%%%%%%%%%%%%%%%%%%%%%%%%%%%%%%%%%%%%%%%%%%%%%%%%%%

Now, using the results prepared in the previous section, we construct functional relations for $\zeta_2({\boldsymbol{s}};G_2)$ and $\zeta(s)$. First we recall the relation for zeta-functions of $C_2$-type which was proved in \cite{KMTJC}, and will extend this relation to that of $G_2$-type. The technique is essentially introduced in our previous papers (see \cite[Remark 7.5]{KMTJC}). 
From \cite[(8.4)]{KMTJC}, we have
\begin{align*}
& \sum_{l\in \mathbb{Z},\,l\not=0 \atop {m\geq 1 \atop {l+m\not=0 \atop l+2m\not=0}}} \frac{(-1)^{l}x^m e^{i(l+2m)\theta}}{l^{2p}m^s(l+m)^{2q}(l+2m)^{2r}}\\
& \ -2\sum_{j=0}^{p}\ \phi(2p-2j)\ \sum_{\xi=0}^{2j} \sum_{\omega=0}^{2j-\xi} \binom{\omega+2r-1}{\omega} (-1)^\omega  \\
& \hspace{0.5in} \times \binom{2q-1+2j-\xi-\omega}{2q-1}(-1)^{2j-\xi-\omega}\frac{1}{2^{2r+\omega}} \sum_{m=1}^\infty \frac{x^m e^{2im\theta}}{m^{s+2q+2j-\xi+2r}}\frac{(i\theta)^{\xi}}{\xi!}  \\
& \ \ -2\sum_{j=0}^{q}\ \phi(2q-2j)\ \sum_{\xi=0}^{2j} \sum_{\omega=0}^{2j-\xi} \binom{\omega+2r-1}{\omega} (-1)^\omega  \\
& \hspace{0.5in} \times \binom{2p-1+2j-\xi-\omega}{2p-1}\sum_{m=1}^\infty \frac{(-1)^m x^m e^{im\theta}}{m^{s+2p+2j-\xi+2r}}\frac{(i\theta)^{\xi}}{\xi!}  \\
& \ \ +2\sum_{j=0}^{r}\ \phi(2r-2j)\ \sum_{\xi=0}^{2j} \sum_{\omega=0}^{2p-1} \binom{\omega+2j-\xi}{\omega} (-1)^\omega  \\
& \hspace{0.5in} \times \binom{2p+2q-2-\omega}{2q-1}(-1)^{2p-1-\omega}\frac{1}{2^{2j-\xi+\omega+1}} \sum_{m=1}^\infty \frac{x^m}{m^{s+2q+2j-\xi+2p}}\frac{(i\theta)^{\xi}}{\xi!}  \\
& \ \ +2\sum_{j=0}^{r}\ \phi(2r-2j)\ \sum_{\xi=0}^{2j} \sum_{\omega=0}^{2q-1} \binom{\omega+2j-\xi}{\omega} (-1)^\omega  \\
& \hspace{0.5in} \times \binom{2p+2q-2-\omega}{2p-1} \sum_{m=1}^\infty \frac{x^m}{m^{s+2p+2j-\xi+2q}}\frac{(i\theta)^{\xi}}{\xi!}=0 
\end{align*}
for $\theta \in [-\pi,\pi]$, $p,q,r \in \mathbb{N}$, $s\in \mathbb{R}$ with $s>1$ and $x\in \mathbb{C}$ with $|x|\leq 1$. Here we use the same method as introduced in our previous papers \cite{KMT,KMTJC} by making use of polylogarithms as follows. Replacing $x$ by $-xe^{i\theta}$ and moving the terms corresponding to $l+3m=0$ of the first member on the left-hand side of the above equation to the right-hand side, we have
\begin{align*}
& \sum_{l\in \mathbb{Z},\,l\not=0 \atop {m \geq 1 \atop {l+m\not=0 \atop {l+2m \not= 0 \atop l+3m\not= 0}}}} \frac{(-1)^{l+m}x^m e^{i(l+3m)\theta}}{l^{2p}m^s(l+m)^{2q}(l+2m)^{2r}}\\
& \ \ -2\sum_{j=0}^{p}\ \phi(2p-2j)\ \sum_{\xi=0}^{2j} \sum_{\omega=0}^{2j-\xi} \binom{\omega+2r-1}{\omega} (-1)^\omega  \\
& \hspace{0.5in} \times \binom{2q-1+2j-\xi-\omega}{2q-1}(-1)^{2j-\xi-\omega}\frac{1}{2^{2r+\omega}} \sum_{m=1}^\infty \frac{(-1)^m x^m e^{3im\theta}}{m^{s+2q+2j-\xi+2r}}\frac{(i\theta)^{\xi}}{\xi!}  \\
& \ \ -2\sum_{j=0}^{q}\ \phi(2q-2j)\ \sum_{\xi=0}^{2j} \sum_{\omega=0}^{2j-\xi} \binom{\omega+2r-1}{\omega} (-1)^\omega  \\
& \hspace{0.5in} \times \binom{2p-1+2j-\xi-\omega}{2p-1}\sum_{m=1}^\infty \frac{x^m e^{2im\theta}}{m^{s+2p+2j-\xi+2r}}\frac{(i\theta)^{\xi}}{\xi!}  \\
& \ \ +2\sum_{j=0}^{r}\ \phi(2r-2j)\ \sum_{\xi=0}^{2j} \sum_{\omega=0}^{2p-1} \binom{\omega+2j-\xi}{\omega} (-1)^\omega  \\
& \hspace{0.5in} \times \binom{2p+2q-2-\omega}{2q-1}(-1)^{2p-1-\omega}\frac{1}{2^{2j-\xi+\omega+1}} \sum_{m=1}^\infty \frac{(-1)^m x^m e^{im\theta}}{m^{s+2q+2j-\xi+2p}}\frac{(i\theta)^{\xi}}{\xi!}  \\
& \ \ +2\sum_{j=0}^{r}\ \phi(2r-2j)\ \sum_{\xi=0}^{2j} \sum_{\omega=0}^{2q-1} \binom{\omega+2j-\xi}{\omega} (-1)^\omega  \\
& \hspace{0.5in} \times \binom{2p+2q-2-\omega}{2p-1} \sum_{m=1}^\infty \frac{(-1)^m x^m e^{im\theta}}{m^{s+2p+2j-\xi+2q}}\frac{(i\theta)^{\xi}}{\xi!} \\
& =- \sum_{m =1}^\infty \frac{(-1)^{m}x^m }{3^{2p}2^{2q}m^{2p+2q+2r+s}}.
\end{align*}
If we fix $p,q,r \in \mathbb{N}$ and $s\in \mathbb{R}$ with $s>1$, then we can apply Lemma \ref{L-4-3} to the above equation with $d=2u$. Consequently we have
\begin{equation}
\begin{split}
& \sum_{l\in \mathbb{Z},\, l\not=0 \atop{m \geq 1 \atop {l+m\not=0 \atop {l+2m \not= 0 \atop l+3m\not= 0}}}} \frac{(-1)^{l+m}x^m e^{i(l+3m)\theta}}{l^{2p}m^s(l+m)^{2q}(l+2m)^{2r}(l+3m)^{2u}}\\
& \ \ +J_1(\theta;x)+J_2(\theta;x)+J_3(\theta;x)+J_4(\theta;x)=0,
\end{split}
\label{eq-4-7} 
\end{equation}
where
\begin{align*}
& J_1 (\theta;x) \\
& =-2\sum_{j=0}^{p}\ \phi(2p-2j)\ \sum_{\xi=0}^{2j} \sum_{\rho=0}^{2j-\xi}\binom{\rho+2u-1}{\rho}(-1)^\rho \sum_{\omega=0}^{2j-\xi-\rho} \binom{\omega+2r-1}{\omega} (-1)^\omega  \\
& \hspace{0.2in} \times 3^{-2u-\rho}\binom{2q-1+2j-\xi-\rho-\omega}{2q-1}\frac{(-1)^{2j-\xi-\rho-\omega}}{2^{2r+\omega}} \sum_{m=1}^\infty \frac{(-1)^m x^m e^{3im\theta}}{m^{s+2q+2r+2u+2j-\xi}}\frac{(i\theta)^{\xi}}{\xi!} \\
& +2\sum_{j=0}^{u}\ \phi(2u-2j)\ \sum_{\xi=0}^{2j} \sum_{\rho=0}^{2p-1}\binom{\rho+2j-\xi}{\rho}(-1)^\rho \sum_{\omega=0}^{2p-1-\rho} \binom{\omega+2r-1}{\omega} (-1)^\omega  \\
& \hspace{0.2in} \times 3^{-2j+\xi-\rho-1}\binom{2p+2q-2-\rho-\omega}{2q-1}\frac{(-1)^{2p-1-\rho-\omega}}{2^{2r+\omega}} \sum_{m=1}^\infty \frac{x^m}{m^{s+2p+2q+2r+2j-\xi}}\frac{(i\theta)^{\xi}}{\xi!}.
\end{align*}
We can similarly write $J_2(\theta;x)$, $J_3(\theta;x)$ and $J_4(\theta;x)$, but they are omitted for the purpose of saving space. 

Next, setting $x=\pm i e^{-3i\theta/2}$ in \eqref{eq-4-7} and 
moving the terms corresponding to $2l+3m=0$ of the first member on the left-hand side to the right-hand side, we have
\begin{equation}
\begin{split}
& \sum_{l\in \mathbb{Z},\,l\not=0 \atop{m \geq 1 \atop {l+m\not=0 \atop {l+2m \not= 0 \atop {l+3m\not= 0 \atop 2l+3m \not=0}}}}} \frac{(-1)^{l+m} (\pm i)^{m} e^{i(2l+3m)\theta/2}}{l^{2p}m^s(l+m)^{2q}(l+2m)^{2r}(l+3m)^{2u}}\\
& \ \ +J_1(\theta;\pm i e^{-3i\theta/2})+J_2(\theta;\pm i e^{-3i\theta/2})+J_3(\theta;\pm i e^{-3i\theta/2})+J_4(\theta;\pm i e^{-3i\theta/2}) \\
& \ \ =\sum_{l,m =1 \atop 2l=3m}^\infty \frac{1}{l^{2p}m^s(-l+m)^{2q}(-l+2m)^{2r}(-l+3m)^{2u}}.
\end{split}
\label{eq-4-9} 
\end{equation}
Note that $(-1)^{l+m} (\pm i)^{m}=(\pm i)^{2l+3m}$. 

Now we apply Lemma \ref{L-4-4} to \eqref{eq-4-9} with $d=2v$ for $v \in \mathbb{N}$. 
In fact, we can see that the left-hand side of \eqref{eq-4-9} is of the same form as \eqref{eq-4-10}. Furthermore, by the same method as in \cite[Section 7]{KMTJC}, we can confirm that 
\begin{equation*}
\begin{split}
& \sum_{l,m =1 \atop {l\not=m \atop {l \not= 2m \atop {l\not= 3m \atop 2l \not=3m}}}}^\infty \frac{1}{l^{2p}m^s(-l+m)^{2q}(-l+2m)^{2r}(-l+3m)^{2u}(-2l+3m)^{2v}} \\
& =\zeta_2(2p,2q,s,2r,2v,2u;G_2)+\zeta_2(2u,2r,s,2q,2v,2p;G_2) \\
& \ \ +\zeta_2(2u,s,2r,2q,2p,2v;G_2)+\zeta_2(2v,2r,2q,s,2u,2p;G_2)\\
& \ \ +\zeta_2(2v,2q,2r,s,2p,2u;G_2).
\end{split}
\end{equation*}
From these results and Theorem \ref{T-3-1}, we obtain the following theorem. 

\begin{theorem} \label{Fn-Rel}
For $p,q,r,u,v \in \mathbb{N}$,
\begin{equation}
\begin{split}
& \zeta_2(2p,s,2q,2r,2u,2v;G_2)+\zeta_2(2p,2q,s,2r,2v,2u;G_2)\\
& \ \ +\zeta_2(2u,2r,s,2q,2v,2p;G_2)+\zeta_2(2u,s,2r,2q,2p,2v;G_2) \\
& \ \ +\zeta_2(2v,2r,2q,s,2u,2p;G_2)+\zeta_2(2v,2q,2r,s,2p,2u;G_2) \\
& \ \ +I_1+I_2+\cdots+I_8=0
\end{split}
\label{Fn-G2}
\end{equation}
holds for all $s \in \mathbb{C}$ except for singularities of functions on the left-hand side, where $I_1,I_2,\ldots,I_8$ are, by using the notation $\phi(s)=(2^{1-s}-1)\zeta(s)$, defined as follows:
\begin{align*}
I_1 & =-2\sum_{k=0}^{p}\zeta(2k)\sum_{\sigma=0}^{2p-2k}\binom{\sigma+2v-1}{\sigma}\sum_{\rho =0}^{2p-2k-\sigma}\binom{\rho +2u-1}{\rho } \\
& \quad\quad\times \sum_{\omega =0}^{2p-2k-\sigma-\rho }\binom{\omega +2r-1}{\omega }\binom{2p+2q-1-2k-\sigma-\rho -\omega }{2q-1}\\
& \quad\quad\times {2^{\sigma-2r-\omega }3^{-2u-2v-\sigma-\rho }} \zeta(s+2p+2q+2r+2u+2v-2k)\\
& \ -2\sum_{k=0}^{v}2^{-2k}\zeta(2k)\sum_{\sigma=0}^{2p-1}\binom{\sigma+2v-2k}{\sigma}\sum_{\rho =0}^{2p-1-\sigma}\binom{\rho +2u-1}{\rho } \\
& \quad\quad\times \sum_{\omega =0}^{2p-1-\sigma-\rho }\binom{\omega +2r-1}{\omega } \binom{2p+2q-2-\sigma-\rho -\omega }{2q-1}{2^{\sigma-2r-\omega }3^{-2u-2v+2k-\sigma-\rho -1}} \\
& \quad\quad\times \{ \zeta(s+2p+2q+2r+2u+2v-2k)+\phi(s+2p+2q+2r+2u+2v-2k)\}\\
& \ -2\sum_{k=0}^{v}\left( 1-2^{-2k} \right)\zeta(2k)\sum_{\sigma=0}^{2p-1}\binom{\sigma+2v-2k}{\sigma}\sum_{\rho =0}^{2p-1-\sigma}\binom{\rho +2u-1}{\rho } \\
& \quad\quad\times \sum_{\omega =0}^{2p-1-\sigma-\rho }\binom{\omega +2r-1}{\omega } \binom{2p+2q-2-\sigma-\rho -\omega }{2q-1} {2^{\sigma-2r-\omega }3^{-2u-2v+2k-\sigma-\rho -1}} \\
& \quad\quad\times \{ \zeta(s+2p+2q+2r+2u+2v-2k)-\phi(s+2p+2q+2r+2u+2v-2k)\};
\end{align*}
\begin{align*}
I_2 & =-2\sum_{k=0}^{u}\zeta(2k)\sum_{\sigma=0}^{2u-2k}\binom{\sigma+2v-1}{\sigma}\sum_{\rho =0}^{2p-1}\binom{\rho +2u-2k-\sigma}{\rho } \\
& \quad\quad\times \sum_{\omega =0}^{2p-2k-\rho }\binom{\omega +2r-1}{\omega }\binom{2p+2q-2-\rho -\omega }{2q-1}{2^{\sigma-2r-\omega }3^{-2u-2v+2k-\rho -1}} \\
& \quad\quad\times \zeta(s+2p+2q+2r+2u+2v-2k)\\
& \ -2\sum_{k=0}^{v}2^{-2k}\zeta(2k)\sum_{\sigma=0}^{2u-1}\binom{\sigma+2v-2k}{\sigma}\sum_{\rho =0}^{2p-1}\binom{\rho +2u-1-\sigma}{\rho } \\
& \quad\quad\times \sum_{\omega =0}^{2p-1-\rho }\binom{\omega +2r-1}{\omega } \binom{2p+2q-2-\rho -\omega }{2q-1}{2^{\sigma-2r-\omega }3^{-2u-2v+2k-\rho -1}} \\
& \quad\quad\times \{ \zeta(s+2p+2q+2r+2u+2v-2k)+\phi(s+2p+2q+2r+2u+2v-2k)\}\\
& \ -2\sum_{k=0}^{v}\left( 1-2^{-2k} \right)\zeta(2k)\sum_{\sigma=0}^{2u-1}\binom{\sigma+2v-2k}{\sigma}\sum_{\rho =0}^{2p-1}\binom{\rho +2u-1-\sigma}{\rho } \\
& \quad\quad\times \sum_{\omega =0}^{2p-1-\rho }\binom{\omega +2r-1}{\omega } \binom{2p+2q-2-\rho -\omega }{2q-1}{2^{\sigma-2r-\omega }3^{-2u-2v+2k-\rho -1}} \\
& \quad\quad\times \{ \zeta(s+2p+2q+2r+2u+2v-2k)-\phi(s+2p+2q+2r+2u+2v-2k)\};
\end{align*}
\begin{align*}
I_3 & =-2\sum_{k=0}^{q}\zeta(2k)\sum_{\sigma=0}^{2q-2k}\binom{\sigma+2v-1}{\sigma}\sum_{\rho =0}^{2q-2k-\sigma}\binom{\rho +2u-1}{\rho } \\
& \quad\quad\times \sum_{\omega =0}^{2q-2k-\sigma-\rho }\binom{\omega +2r-1}{\omega }\binom{2p+2q-1-2k-\sigma-\rho -\omega }{2p-1}(-1)^{\sigma+\rho +\omega }{2^{\sigma-2u-\rho }} \\
& \quad\quad\times \zeta(s+2p+2q+2r+2u+2v-2k)\\
& \ +2\sum_{k=0}^{v}2^{-2k}\zeta(2k)\sum_{\sigma=0}^{2q-1}\binom{\sigma+2v-2k}{\sigma}\sum_{\rho =0}^{2q-1-\sigma}\binom{\rho +2u-1}{\rho } \\
& \quad\quad\times \sum_{\omega =0}^{2q-1-\sigma-\rho }\binom{\omega +2r-1}{\omega } \binom{2p+2q-2-\sigma-\rho -\omega }{2p-1}(-1)^{\sigma+\rho +\omega }{2^{\sigma-2u-\rho }} \\
& \quad\quad\times \{ \zeta(s+2p+2q+2r+2u+2v-2k)+\phi(s+2p+2q+2r+2u+2v-2k)\}\\
& \ +2\sum_{k=0}^{v}\left( 1-2^{-2k} \right)\zeta(2k)\sum_{\sigma=0}^{2q-1}\binom{\sigma+2v-2k}{\sigma}\sum_{\rho =0}^{2q-1-\sigma}\binom{\rho +2u-1}{\rho } \\
& \quad\quad\times \sum_{\omega =0}^{2q-1-\sigma-\rho }\binom{\omega +2r-1}{\omega } \binom{2p+2q-2-\sigma-\rho -\omega }{2p-1}(-1)^{\sigma+\rho +\omega }{2^{\sigma-2u-\rho }} \\
& \quad\quad\times \{ \zeta(s+2p+2q+2r+2u+2v-2k)-\phi(s+2p+2q+2r+2u+2v-2k)\};
\end{align*}
\begin{align*}
I_4 & =2\sum_{k=0}^{u}\zeta(2k)\sum_{\sigma=0}^{2u-2k}\binom{\sigma+2v-1}{\sigma}\sum_{\rho =0}^{2q-1}\binom{\rho +2u-2k-\sigma}{\rho } \\
& \quad\quad\times \sum_{\omega =0}^{2q-1-\rho }\binom{\omega +2r-1}{\omega }\binom{2p+2q-2-\rho -\omega }{2q-1}\\
& \quad\quad\times (-1)^{\rho +\omega }{2^{-2u+2k+2\sigma-\rho -1}3^{-2v-\sigma}} \zeta(s+2p+2q+2r+2u+2v-2k)\\
& \ +2\sum_{k=0}^{v}2^{-2k}\zeta(2k)\sum_{\sigma=0}^{2u-1}\binom{\sigma+2v-2k}{\sigma}\sum_{\rho =0}^{2q-1}\binom{\rho +2u-1-\sigma}{\rho } \\
& \quad\quad\times \sum_{\omega =0}^{2q-1-\rho }\binom{\omega +2r-1}{\omega } \binom{2p+2q-2-\rho -\omega }{2p-1}\\
& \quad\quad\times (-1)^{\rho +\omega }{2^{-2u+2\sigma-\rho }3^{-2v+2k-\sigma-1}} \\
& \quad\quad\times \{ \zeta(s+2p+2q+2r+2u+2v-2k)+\phi(s+2p+2q+2r+2u+2v-2k)\}\\
& \ +2\sum_{k=0}^{v}\left( 1-2^{-2k} \right)\zeta(2k)\sum_{\sigma=0}^{2u-1}\binom{\sigma+2v-2k}{\sigma}\sum_{\rho =0}^{2q-1}\binom{\rho +2u-1-\sigma}{\rho } \\
& \quad\quad\times \sum_{\omega =0}^{2q-1-\rho }\binom{\omega +2r-1}{\omega } \binom{2p+2q-2-\rho -\omega }{2p-1}\\
& \quad\quad\times (-1)^{\rho +\omega }{2^{-2u+2\sigma-\rho }3^{-2v+2k-\sigma-1}} \\
& \quad\quad\times \{ \zeta(s+2p+2q+2r+2u+2v-2k)-\phi(s+2p+2q+2r+2u+2v-2k)\};
\end{align*}
\begin{align*}
I_5 & =-2\sum_{k=0}^{r}\zeta(2k)\sum_{\sigma=0}^{2r-2k}\binom{\sigma+2v-1}{\sigma}\sum_{\rho =0}^{2r-2k-\sigma}\binom{\rho +2u-1}{\rho } \\
& \quad\quad\times \sum_{\omega =0}^{2p-1}\binom{\omega +2r-2k-\sigma-\rho }{\omega }\binom{2p+2q-2-\omega }{2q-1}\\
& \quad\quad\times (-1)^\rho {2^{-2r+2k+2\sigma+\rho -\omega -1}} \zeta(s+2p+2q+2r+2u+2v-2k)\\
& \ -2\sum_{k=0}^{v}2^{-2k}\zeta(2k)\sum_{\sigma=0}^{2r-1}\binom{\sigma+2v-2k}{\sigma}\sum_{\rho =0}^{2r-1-\sigma}\binom{\rho +2u-1}{\rho } \\
& \quad\quad\times \sum_{\omega =0}^{2p-1}\binom{\omega +2r-1-\sigma-\rho }{\omega } \binom{2p+2q-2-\omega }{2q-1}(-1)^\rho  {2^{-2r+2\sigma+\rho -\omega }} \\
& \quad\quad\times \{ \zeta(s+2p+2q+2r+2u+2v-2k)+\phi(s+2p+2q+2r+2u+2v-2k)\}\\
& \ -2\sum_{k=0}^{v}\left( 1-2^{-2k} \right)\zeta(2k)\sum_{\sigma=0}^{2r-1}\binom{\sigma+2v-2k}{\sigma}\sum_{\rho =0}^{2p-1-\sigma}\binom{\rho +2u-1}{\rho } \\
& \quad\quad\times \sum_{\omega =0}^{2p-1}\binom{\omega +2r-1-\sigma-\rho}{\omega } \binom{2p+2q-2-\omega }{2q-1}(-1)^\rho {2^{-2r+2\sigma+\rho -\omega }} \\
& \quad\quad\times \{ \zeta(s+2p+2q+2r+2u+2v-2k)-\phi(s+2p+2q+2r+2u+2v-2k)\};
\end{align*}
\begin{align*}
I_6 & =2\sum_{k=0}^{u}\zeta(2k)\sum_{\sigma=0}^{2u-2k}\binom{\sigma+2v-1}{\sigma}\sum_{\rho =0}^{2r-1}\binom{\rho +2u-2k-\sigma}{\rho } \\
& \quad\quad\times \sum_{\omega =0}^{2p-1}\binom{\omega +2r-1-\rho }{\omega }\binom{2p+2q-2-\omega }{2q-1}\\
& \quad\quad\times (-1)^\rho {2^{-2r+\sigma+\rho -\omega }3^{-2v-\sigma}} \zeta(s+2p+2q+2r+2u+2v-2k)\\
& \ +2\sum_{k=0}^{v}2^{-2k}\zeta(2k)\sum_{\sigma=0}^{2u-1}\binom{\sigma+2v-2k}{\sigma}\sum_{\rho =0}^{2r-1}\binom{\rho +2u-1-\sigma}{\rho } \\
& \quad\quad\times \sum_{\omega =0}^{2p-1}\binom{\omega +2r-1-\rho }{\omega } \binom{2p+2q-2-\omega }{2q-1}\\
& \quad\quad\times (-1)^\rho {2^{-2r+\sigma+\rho -\omega }3^{-2v-\sigma}} \\
& \quad\quad\times \{ \zeta(s+2p+2q+2r+2u+2v-2k)+\phi(s+2p+2q+2r+2u+2v-2k)\}\\
& \ +2\sum_{k=0}^{v}\left( 1-2^{-2k} \right)\zeta(2k)\sum_{\sigma=0}^{2u-1}\binom{\sigma+2v-2k}{\sigma}\sum_{\rho =0}^{2r-1}\binom{\rho +2u-1-\sigma}{\rho } \\
& \quad\quad\times \sum_{\omega =0}^{2p-1}\binom{\omega +2r-1-\rho }{\omega } \binom{2p+2q-2-\omega }{2q-1}\\
& \quad\quad\times (-1)^\rho {2^{-2r+\sigma+\rho -\omega }3^{-2v-\sigma}} \\
& \quad\quad\times \{ \zeta(s+2p+2q+2r+2u+2v-2k)-\phi(s+2p+2q+2r+2u+2v-2k)\};
\end{align*}
\begin{align*}
I_7 & =2\sum_{k=0}^{r}\zeta(2k)\sum_{\sigma=0}^{2r-2k}\binom{\sigma+2v-1}{\sigma}\sum_{\rho =0}^{2r-2k-\sigma}\binom{\rho +2u-1}{\rho } \\
& \quad\quad\times \sum_{\omega =0}^{2q-1}\binom{\omega +2r-2k-\sigma-\rho }{\omega }\binom{2p+2q-2-\omega }{2p-1}(-1)^{\rho +\omega }{2^{\sigma}} \\
& \quad\quad\times \zeta(s+2p+2q+2r+2u+2v-2k)\\
& \ +2\sum_{k=0}^{v}2^{-2k}\zeta(2k)\sum_{\sigma=0}^{2r-1}\binom{\sigma+2v-2k}{\sigma}\sum_{\rho =0}^{2r-1-\sigma}\binom{\rho +2u-1}{\rho } \\
& \quad\quad\times \sum_{\omega =0}^{2q-1}\binom{\omega +2r-1-\sigma-\rho }{\omega } \binom{2p+2q-2-\omega }{2p-1}(-1)^{\rho +\omega }{2^{\sigma}} \\
& \quad\quad\times \{ \zeta(s+2p+2q+2r+2u+2v-2k)+\phi(s+2p+2q+2r+2u+2v-2k)\}\\
& \ +2\sum_{k=0}^{v}\left( 1-2^{-2k} \right)\zeta(2k)\sum_{\sigma=0}^{2r-1}\binom{\sigma+2v-2k}{\sigma}\sum_{\rho =0}^{2r-1-\sigma}\binom{\rho +2u-1}{\rho } \\
& \quad\quad\times \sum_{\omega =0}^{2q-1}\binom{\omega +2r-1-\sigma-\rho }{\omega } \binom{2p+2q-2-\omega }{2p-1}(-1)^{\rho +\omega }{2^{\sigma}} \\
& \quad\quad\times \{ \zeta(s+2p+2q+2r+2u+2v-2k)-\phi(s+2p+2q+2r+2u+2v-2k)\};
\end{align*}
\begin{align*}
I_8 & =-2\sum_{k=0}^{u}\zeta(2k)\sum_{\sigma=0}^{2u-2k}\binom{\sigma+2v-1}{\sigma}\sum_{\rho =0}^{2r-1}\binom{\rho +2u-2k-\sigma}{\rho } \\
& \quad\quad\times \sum_{\omega =0}^{2q-1}\binom{\omega +2r-1-\rho }{\omega }\binom{2p+2q-2-\omega }{2p-1}\\
& \quad\quad\times (-1)^{\rho +\omega }{2^{\sigma}3^{-2v-\sigma}} \zeta(s+2p+2q+2r+2u+2v-2k)\\
& \ -2\sum_{k=0}^{v}2^{-2k}\zeta(2k)\sum_{\sigma=0}^{2u-1}\binom{\sigma+2v-2k}{\sigma}\sum_{\rho =0}^{2r-1}\binom{\rho +2u-1-\sigma}{\rho } \\
& \quad\quad\times \sum_{\omega =0}^{2q-1}\binom{\omega +2r-1-\rho }{\omega } \binom{2p+2q-2-\omega }{2p-1}\\
& \quad\quad\times (-1)^{\rho +\omega }{2^{\sigma}3^{-2v+2k-\sigma-1}} \\
& \quad\quad\times \{ \zeta(s+2p+2q+2r+2u+2v-2k)+\phi(s+2p+2q+2r+2u+2v-2k)\}\\
& \ -2\sum_{k=0}^{v}\left( 1-2^{-2k} \right)\zeta(2k)\sum_{\sigma=0}^{2u-1}\binom{\sigma+2v-2k}{\sigma}\sum_{\rho =0}^{2r-1}\binom{\rho +2u-1-\sigma}{\rho } \\
& \quad\quad\times \sum_{\omega =0}^{2q-1}\binom{\omega +2r-1-\rho}{\omega } \binom{2p+2q-2-\omega }{2p-1}\\
& \quad\quad\times (-1)^{\rho +\omega }{2^{\sigma}3^{-2v+2k-\sigma-1}} \\
& \quad\quad\times \{ \zeta(s+2p+2q+2r+2u+2v-2k)-\phi(s+2p+2q+2r+2u+2v-2k)\}.
\end{align*}
\end{theorem}

\begin{example} Putting $(p,q,r,u,v)=(1,1,1,1,1)$ in \eqref{Fn-G2}, we have
\begin{equation}
\begin{split}
& \zeta_{2}(2,s,2,2,2,2;G_2)+\zeta_{2}(2,2,s,2,2,2;G_2)+\zeta_{2}(2,2,2,s,2,2;G_2) \\
& \ \ = -\frac{5}{1458}\left( 2^{-s}+\frac{5519}{4}\right) \zeta(s+10) -\frac{1}{162} \left( 2^{-s}-466\right) \zeta(2)\zeta(s+8).
\end{split}
\label{Eq-G_2}
\end{equation}
In particular when $s=2$, we recover 
\begin{align*}
&  \zeta_{2}(2,2,2,2,2,2;G_2) = \frac{23}{297904566960} \pi^{12}, 
\end{align*}
which was already obtained in Example \ref{Exam-2-1}. 
\end{example}

\begin{remark} \label{Rem-4-7}
In \cite{Zhao}, Zhao expressed several values $\zeta_2({\bf k};G_2)$ for ${\bf k}\in \mathbb{N}_0^6$ in terms of double polylogarithms and gave approximate values of them, for example, 
\begin{align*}
\zeta_2(2,1,1,1,1,1;G_2) & =  0.0099527234\cdots.
\end{align*}
By using the same method as stated above, we can explicitly obtain
\begin{equation}
\zeta_2(2,1,1,1,1,1;G_2) = -\frac{109}{1296}\zeta(7)+\frac{1}{18}\zeta(2)\zeta(5), \label{Zhao-1}
\end{equation}
which agrees with Zhao's numerical computation. 
We can further give a functional relation between $\zeta_2({\bf s};G_2)$ and $\zeta(s)$ including \eqref{Zhao-1}, which is an analogue of \eqref{Eq-G_2}. 
Under more preparations, we will be able to give evaluation formulas for a certain class of values $\zeta_2({\bf k};G_2)$ for ${\bf k}\in \mathbb{N}_0^6$ in terms of $\zeta(m)$ for $m \geq 2$. We will state the details in a forthcoming paper.
\end{remark}

\bigskip

\baselineskip 14pt

\bibliographystyle{amsplain}

\ 

\end{document}